\documentclass[12pt]{article}
\usepackage{amssymb,latexsym}
\input epsf

\title{\textbf{The absolute order of a permutation representation
of a Coxeter group}}

\author{Christos~A.~Athanasiadis\\
Department of Mathematics\\
University of Athens\\
Athens 15784, Hellas (Greece)\\
Email: \texttt{caath@math.uoa.gr}
\vspace{5pt}
\and
Yuval~Roichman\\
Department of Mathematics
and Computer Science\\
Bar-Ilan University\\
Ramat Gan\\
52900, Israel\\
Email: \texttt{yuvalr@math.biu.ac.il} }

\date{\small February 27, 2013}

  \def\CC{{\mathbb C}}

  \def\NN{{\mathbb N}}

  \def\ZZ{{\mathbb Z}}
  \def\RR{{\mathbb R}}
  
  \def\aA{{\mathcal A}}
  
  \def\hH{{\mathcal H}}
  
  \def\lL{{\mathcal L}}
  \def\mM{{\mathcal M}}
  
  \def\Abs{{\rm Abs}}

  \def\Fix{{\rm Fix}}
  \def\Mov{{\rm Mov}}
  \def\rank{{\rm rank}}
  \def\rk{{\rm rk}}

  \def\s{{\tt{S}}}
  \def\T{{\tt{T}}}
  \newcommand{\qed}{$\hfill \Box$}

\begin{document}
\maketitle

\newtheorem{theorem}{Theorem}[section]
\newtheorem{proposition}[theorem]{Proposition}
\newtheorem{corollary}[theorem]{Corollary}
\newtheorem{defn}[theorem]{Definition}
\newtheorem{remark}[theorem]{Remark}
\newtheorem{lemma}[theorem]{Lemma}
\newtheorem{example}[theorem]{Example}
\newtheorem{examples}[theorem]{Examples}
\newtheorem{conjecture}[theorem]{Conjecture}
\newtheorem{fact}[theorem]{Fact}
\newtheorem{question}[theorem]{Question}
\newtheorem{observation}[theorem]{Observation}
\newtheorem{claim}[theorem]{Claim}

\begin{abstract}
A permutation representation of a Coxeter group $W$ naturally
defines an absolute order. This family of partial orders (which
includes the absolute order on $W$) is introduced and studied in
this paper. Conditions under which the associated rank generating
polynomial divides the rank generating polynomial of the absolute
order on $W$ are investigated when $W$ is finite. Several
examples, including a symmetric group action on perfect matchings,
are discussed. As an application, a well-behaved absolute order
on the alternating subgroup of $W$ is defined.

\end{abstract}


\section{Introduction}

The Bruhat order on a Coxeter group $W$ is a key ingredient in
understanding the structure of $W$. This order involves both the
set of simple reflections $\s$ and the set of all reflections $\T$
of $W$: it may be defined by the condition that $u \in W$ is covered
by $v \in W$ if there exists $t \in \T$ such that $v = tu$ and
$\ell_\s (v) = \ell_\s (u) + 1$, where $\ell_\s: W \to \NN$ is
the length function with respect to the generating set $\s$. There
are two ``more coherent" closely related concepts. Replacing the
role of $\T$ by $\s$ determines an order which was extensively
studied in the past three decades, namely the weak order on $W$.
Replacing the role of $\s$ by $\T$ determines the absolute order.
The study of maximal chains in the absolute order on the symmetric
group is traced at least back to Hurwitz~\cite{Hurwitz}; see
also~\cite{Denes, Strehl}. However, the growing interest in the
absolute order is relatively recent and followed the discovery
\cite{Be, BW} that distinguished intervals in the absolute order,
known as the noncrossing partition lattices, are objects of
importance in the theory of finite-type Artin groups. For further
information on the absolute order, the reader is referred
to~\cite[Section 2.4]{Armstrong} \cite{AK, Ka}.

Consider a transitive action of $W$ on a set $X$. Motivated by
recent work of Rains and Vazirani~\cite{Rains-Vazirani}, which
introduces and studies the Bruhat order on $X$, a naturally defined
absolute order on $X$ is introduced in this paper. Our goal is to
find conditions under which important enumerative and structural
properties of the absolute order on the acting group $W$ carry over
to the absolute order on $X$; in particular, conditions under which
the associated rank generating polynomial divides the rank generating
polynomial of the absolute order on $W$. Several examples, including
the symmetric group action on ordered tuples and its conjugation
action on fixed point free involutions, are discussed. As an
application, a well-behaved absolute order on the alternating
subgroup of $W$ is defined and studied.

\section{Basic concepts}
\label{sec:basic}

Let $W$ be a Coxeter group with set of reflections $\T$ (for
background on Coxeter groups the reader is referred to \cite{BjB,
Bou, Humphreys}). The minimum length of a $\T$-word for an element
$w \in W$ is denoted by $\ell_\T (w)$ and called the absolute
length of $w$. The absolute order on $W$, denoted by $\Abs(W)$, is
the partial order $(W, \le_\T)$ defined by letting $u \le_\T v$ if
$\ell_\T (vu^{-1}) = \ell_\T (v) - \ell_\T (u)$, for $u, v \in W$.
Equivalently, $\le_\T$ is the reflexive and transitive closure of
the relation on $W$ consisting of the pairs $(u, v)$ of elements
of $W$ for which $\ell_\T (u) < \ell_\T (v)$ and $v = tu$ for some
$t \in \T$. For basic properties of $\Abs(W)$, see~\cite[Section
2.4]{Armstrong}.

We will be concerned with the following generalization of the absolute
order on $W$. Consider a transitive action $\rho$ of $W$ on a set $X$.
We will write $wx$ for the result $\rho(w) (x)$ of the action of $w \in
W$ on an element $x \in X$.

\begin{defn} \label{def:main1} \rm
Fix an arbitrary element $x_0 \in X$.
\begin{itemize}
\itemsep=0pt
\item[(a)]
The absolute length of $x \in X$ is defined as $\ell_\T (x) := \min \,
\{\ell_\T (w): x = w x_0\}$.

\item[(b)]
The absolute order on $X$, denoted $\Abs(X)$, associated to $\rho$ is
the partial order $(X, \le_\T)$ defined by letting $x \le_\T y$ if there
exists $w \in W$ such that $y = wx$ and $\ell_\T (w) = \ell_\T (y) -
\ell_\T (x)$, for $x, y \in X$. Equivalently, $\le_\T$ is the reflexive
and transitive closure of the relation on $X$ consisting of the pairs
$(x, y)$ of elements of $X$ for which $\ell_\T (x) < \ell_\T (y)$ and
$y = tx$ for some $t \in \T$.
\end{itemize}
\end{defn}

The present section discusses elementary properties and examples of
$\Abs(X)$. We begin with some comments on Definition \ref{def:main1}.

\begin{remark} \label{rem:main} \rm
(a) A different way to describe the relation $\le_\T$ on $X$ is
the following. Let $x_0 \in X$ be fixed, as before, and consider
the simple graph $\Gamma = \Gamma(W,\rho)$ on the vertex set $X$
whose (undirected) edges are the sets of the form $\{x, tx\}$ for
$t \in \T$ and $x \in X$. Then for every $x \in X$, the absolute
length $\ell_\T (x)$ is equal to the distance between $x_0$ and
$x$ in the graph $\Gamma$ and for $x, y \in X$, we have $x \le_\T
y$ if and only if $x$ lies in a geodesic path in $\Gamma$ with
endpoints $x_0$ and $y$. This description implies that $\le_\T$ is
indeed a partial order on $X$ and that it coincides with the
reflexive and transitive closure of the relation on $X$ described
in Definition \ref{def:main1} (b).

(b) The isomorphism type of $\Abs(X)$ is independent of the choice
of $x_0 \in X$. Indeed, consider another base point $y_0 \in X$
and let $\Abs(X, x_0)$ and $\Abs(X, y_0)$ denote the absolute
orders on $X$ with respect to $x_0$ and $y_0$, respectively.
Choose $w_0 \in W$ so that $y_0 = w_0 x_0$ and define a map $f: X
\mapsto X$ by letting $f(x) = w_0 x$ for $x \in X$. Clearly, $f$
is a bijection and satisfies $f(x_0) = y_0$. Moreover, since $\T$
is closed under conjugation, the map $f$ is an automorphism of the
graph $\Gamma$ considered in part (a). These properties imply that
$f: \Abs(X, x_0) \mapsto \Abs (X, y_0)$ is an isomorphism of
partially ordered sets.

(c) The order $\Abs(X)$ has minimum element $x_0$.

(d) As an easy consequence of the definition of absolute length, we
have $\ell_\T (wx) \le \ell_\T (w) + \ell_\T (x)$ for all $w \in
W$ and $x \in X$.
 \qed
\end{remark}

Since the action $\rho$ is transitive, the set $X$ may be identified
with the set of left cosets of the stabilizer of $x_0 \in X$ in $W$.
This identification leads to the following reformulation of Definition
\ref{def:main1}, which we will often find convenient (the role of the
base point $x_0$ in Definition \ref{def:main1} will be played by the
subgroup $H$).

\begin{defn}\label{def:main2} \rm
Let $H$ be a subgroup of $W$ and let $X = W/H$ be the set of left cosets
of $H$ in $W$.
\begin{itemize}
\itemsep=0pt
\item[(a)]
The absolute length of $x \in X$ is defined as $\ell_\T (x) := \min \
\{\ell_\T (w): w\in x\}$.

\item[(b)]
The absolute order on $X$, denoted $\Abs(X)$, is the partial order $(X,
\le_\T)$ defined by letting $x \le_\T y$ if there exists $w \in W$ such
that $y = wx$ and $\ell_\T (w) = \ell_\T (y) - \ell_\T (x)$, for $x, y
\in X$.
\end{itemize}
\end{defn}

We recall that a partially ordered set (poset) $P$ with a minimum element
$\hat{0}$ is said to be locally graded with rank function $\rk: P \to \NN$
if for each $x \in P$, every maximal chain in the closed interval
$[\hat{0}, x]$ of $P$ has exactly $\rk (x) + 1$ elements (for background
and terminology on posets we refer to \cite[Chapter 3]{StaEC1}). We
note the following elementary property of $\Abs(X)$.

\begin{proposition} \label{prop:elementary}
The absolute order $\Abs(X)$ is locally graded, with minimum element
$\hat{0} = x_0$ and rank function given by the absolute length.
\end{proposition}

\noindent
\emph{Proof.}
We have already noted that $x_0$ is the minimum element of $\Abs(X)$.
Thus, it suffices to show that $\ell_\T (y) = \ell_\T (x) + 1$
whenever $y$ covers $x$ in $\Abs(X)$. This is an easy consequence of
Definition \ref{def:main1}.
\qed

\bigskip
We recall (see, for instance, \cite[Theorem 2.7.3]{Armstrong}\cite[Section
3.9]{Humphreys} and the references given there) that when $W$ is finite,
the rank (or length) generating polynomial of $\Abs(W)$ satisfies
  \begin{equation}\label{eq:rankgenW}
    W_\T(q) \ := \ \sum\limits_{w \in W} q^{\ell_\T (w)} \ = \
    \prod_{i=1}^d (1 + e_i q),
  \end{equation}
where $d$ is the Coxeter rank of $W$ and $e_1, e_2,\dots,e_d$ are
its exponents. The rank generating polynomial
  \begin{equation}\label{eq:rankgen}
    X_\T(q) \ := \ \sum\limits_{x \in X} q^{\ell_\T (x)}
  \end{equation}
of $\Abs(X)$ is well-defined when $X$ is finite. The following question
provided much of the motivation for this paper.

\begin{question} \label{que:fact}
For which $W$-actions $\rho$ does $X_\T(q)$ divide $W_\T(q)$?
\end{question}

  \begin{figure}
  \epsfysize = 2.0 in
  \centerline{\epsffile{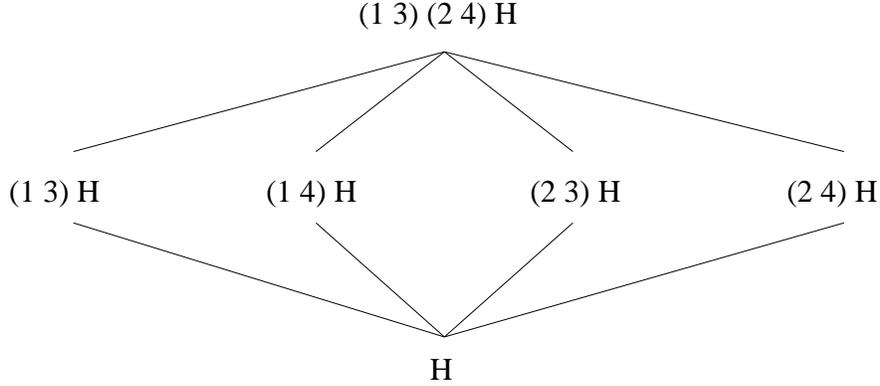}}
  \caption{An absolute order of $S_4$.}
  \label{fig:absX}
  \end{figure}

\medskip
We now list examples, some of which will be studied in detail in later
sections.

\begin{example} \label{ex:abs(W)} \rm
(a) The order $\Abs(W)$ occurs by letting $\rho$ be the left
multiplication action of $W$ on itself and choosing $x_0$ as the
identity element $e \in W$ in Definition \ref{def:main1}, or by choosing
$H$ as the trivial subgroup $\{e\}$ of $W$ in Definition \ref{def:main2}.

(b) Let $H$ be the subgroup of $W$ generated by a given reflection $t_0
\in \T$. The set $X = W/H$ of left cosets of $H$ in $W$ is in bijection
with the alternating subgroup $W^+$ of $W$ and hence $\Abs(X)$ gives
rise to an absolute order on $W^+$. This order will be studied in Section
\ref{sec:alt}.

(c) Let $\lambda$ be an integer partition of $m$ and let $X$ consist of
the set partitions of $\{1, 2,\dots,m\}$ whose block sizes are the parts
of $\lambda$. The symmetric group $S_m$ acts transitively on $X$ and
thus defines an absolute order. This order will be studied in Section
\ref{subsection:PM} in the motivating special case in which $m = 2n$ is
even and all parts of $\lambda$ are equal to 2. The resulting absolute
order is a partial order on the set of perfect matchings of $\{1,
2,\dots,2n\}$. The stabilizer of this action is the natural embedding
of the hyperoctahedral group $B_n$ in $S_{2n}$.

(d) Let $\lambda = (\lambda_1, \lambda_2,\dots,\lambda_r)$ be an
integer partition of $n$ and let $X$ consist of the ordered set
partitions (meaning, set partitions in which the order of the
blocks matters) of $\{1, 2,\dots,n\}$ whose block sizes are
$\lambda_1, \lambda_2,\dots,\lambda_r$, in this order. The
symmetric group $S_n$ acts transitively on $X$ and the stabilizer
is a Young subgroup $S_{\lambda_1} \times \cdots \times
S_{\lambda_r}$ of $S_n$. The resulting absolute order will be
discussed in Section \ref{sec:tuples} in the special case in
which $\lambda = (n-k, 1,\dots,1)$, where $k \in \{1,
2,\dots,n-1\}$. Then $X$ can be identified with the set of
$k$-tuples of pairwise distinct elements of $\{1, 2,\dots,n\}$.

(e) Consider the special case $n = 4$, $\lambda = (2, 2)$ and $x_0
= (\{1, 2\}, \{3, 4\})$ of the example of part (d). Equivalently,
let $W$ be the symmetric group $S_4$ and let $H$ be the four element
subgroup generated by the commuting reflections $(1 \ 2)$ and $(3 \
4)$. Then $X = W/H$ has six elements. The Hasse diagram of $\Abs(X)$
is shown on Figure \ref{fig:absX}.
\qed
\end{example}

\begin{remark} \label{rem:cycles} \rm
It is possible that not all edges of the graph $\Gamma = \Gamma
(W, \rho)$, defined in Remark \ref{rem:main} (a), are edges of the
Hasse diagram of $\Abs(X)$. For instance, consider the action of
$S_4$ on the set $X$ of perfect matchings of $\{1, 2, 3, 4\}$,
discussed in Example \ref{ex:abs(W)} (c). Then $X$ has three
elements and $\Gamma$ is the complete graph on these three
vertices. On the other hand, $\Abs(X)$ has a minimum element $x_0$
which is covered by the other two elements of $X$. Thus exactly
one of the edges of $\Gamma$ is not an edge of the Hasse diagram
of $\Abs(X)$. \qed
\end{remark}

\section{Modular subgroups}
\label{sec:modular}

This section investigates a natural condition on a subgroup of a
Coxeter group, called modularity, and shows that under this condition,
the corresponding absolute order is well-behaved in several ways.
Enumerative (Proposition \ref{prop:fact}) and order-theoretic
(Theorem \ref{thm:modular}) characterizations, as well as examples,
of modularity are given. Throughout this section, $W$ is a Coxeter
group with identity element $e$, $\T$ is the set of reflections, $H$
is a subgroup of $W$ and $X = W/H$ is the set of left cosets of $H$
in $W$. The Coxeter rank of $W$ will be denoted by $\rank(W)$.

The following elementary properties of absolute length (proofs are
left to the reader) will be frequently used throughout this paper.

\begin{fact} \label{absolute-length-properties}
For $u, v, w \in W$ we have:
\begin{itemize}
\itemsep=0pt
\item[{\rm (a)}] $\ell_\T (w) = 0 \Leftrightarrow w=e$,

\item[{\rm (b)}] $\ell_\T (w) = 1 \Leftrightarrow w \in \T$,

\item[{\rm (c)}] $\ell_\T (w^{-1}) = \ell_\T (w)$,

\item[{\rm (d)}] $\ell_\T (uv) \le \ell_\T (u) + \ell_\T (v)$,

\item[{\rm (e)}] $\ell_\T (wuw^{-1}) = \ell_\T (u)$,

\item[{\rm (f)}] $\ell_\T (w) = {\rm codim} \, (\Fix(w))$, if $W$ is
finite,
\end{itemize}
where $\Fix(w)$ is the fixed space of $w$ when $W$ is realized as
a group generated by reflections in Euclidean space (see the relevant
discussions after Remark~\ref{rem:induced}).
\end{fact}

\smallskip

The main definition of this section is as follows.

\begin{defn} \label{def:modular} \rm
We say that $H$ is a \textit{modular subgroup} of $W$ if every left
coset of $H$ in $W$ has a minimum in $\Abs(W)$.
\end{defn}

We note that for $x \in X$ and $w_\circ \in x$, the element $w_\circ$
is the minimum of $x$ in $\Abs(W)$ if and only if we have $\ell_\T
(w_\circ h) = \ell_\T (w_\circ) + \ell_\T (h)$ for every $h \in H$.
We also note that if $H$ is a modular subgroup of $W$, then so are its
conjugate subgroups.

\begin{example} \label{ex:modular1} \rm
(a) Let $H$ be a subgroup of $W$ generated by a single reflection $t
\in \T$. Then every left coset $x \in X$ consists of two elements $w$
and $wt$, which are comparable in $\Abs(W)$. This implies that $H$ is
a modular subgroup of $W$.

(b) Let $H$ be the symmetric group $S_{n-1}$, naturally embedded in
$S_n$. It will be shown in Example \ref{ex:SkBk} (and can be verified
directly) that $H$ is a modular subgroup of $S_n$. The corresponding
absolute order consists of the minimum element $H$ and the left cosets
$(i \ n) \, H$ for $i \in \{1, 2,\dots,n-1\}$, each of which covers
$H$.

(c) The subgroup $H$ of $S_4$ in part (e) of Example \ref{ex:abs(W)}
is not modular. Indeed, there is a single left coset $w H
\in X$, that with $w = (1 \ 3) (2 \ 4)$, which does not have a minimum
in $\Abs(S_4)$. As an induced subposet of $\Abs(W)$, this coset has
$w$ and $(1 \ 4) (2 \ 3)$ as minimal elements, $(1 \ 4 \ 2 \ 3)$ and
$(1 \ 3 \ 2 \ 4)$ as maximal elements and all four possibe Hasse edges
among these elements.

(d) It is possible for a subgroup $H$ of a finite Coxeter group $W$
to have a left coset which has a unique element of minimum absolute
length but no minimum in $\Abs(W)$ (clearly, such a subgroup $H$ cannot
be modular). Consider, for instance, the hyperoctahedral group $B_n$
for some $n \ge 4$ and write $( \! (a \ b) \! )$ for the reflection
in $W$ with cycle form $(a \ b) (-a \, -b)$. Let $H$ be the subgroup
of order 16 generated by the pairwise commuting reflections $t_1 = (
\! (1 \ 2) \! )$, $t_2 = ( \! (1 \, -2) \! )$, $t_3 = ( \! (3 \ 4) \!
)$ and $t_4 = ( \! (3 \, -4) \! )$ and let $t = ( \! (1 \ 3) \! )$
and $h = t_1t_2t_3t_4 \in H$. Then $t H$ contains a unique reflection,
namely $t$, but has no minimum element in $\Abs(W)$, since $t$ is not
comparable to $t h$.
\qed
\end{example}

The following proposition explains the significance of modularity
with respect to Question \ref{que:fact}. It should be compared to
\cite[Lemma~7.1.2]{BjB} \cite[Section~5.2]{Humphreys}
\cite[Theorem~8.1]{Rains-Vazirani}.

\begin{proposition} \label{prop:fact}
Assume that $W$ is finite. Then the subgroup $H$ is modular
if and only if $W_\T(q) = H_\T(q) \cdot X_\T (q)$.
\end{proposition}

\noindent
\emph{Proof.}
Let $w_x \in x$ be an element of minimum absolute length in $x \in
X$. Thus, we have $\ell_\T (w_x) = \ell_\T (x)$ for every $x \in
X$ and hence $\ell_\T (w_x h) \le \ell_\T (w_x) + \ell_\T (h) =
\ell_\T (x) + \ell_\T (h)$ for all $x \in X$ and $h \in H$. As a
result, we find that

  \begin{eqnarray*}
  W_\T(q) &=& \sum\limits_{w \in W} \, q^{\ell_\T(w)} \ = \
  \sum\limits_{x \in X} \sum\limits_{h \in H} \,
  q^{\ell_\T (w_x h)} \\
  & & \\
  & \preceq & \sum\limits_{x \in X} \sum\limits_{h \in H} \,
  q^{\ell_\T (x) \, + \, \ell_\T (h)} \ = \ X_\T (q) \cdot H_\T(q),
  \end{eqnarray*}
where $\preceq$ stands for the reverse lexicographic order on the
set of polynomials with nonnegative integer coefficients, i.e.,
for $f(q), g(q) \in \NN[q]$ we write $f(q) \prec g(q)$ if the
highest term of $g(q) - f(q)$ has positive coefficient. Equality
holds if and only if $\ell_\T (w_x h) = \ell_\T (w_x) + \ell_\T (h)$,
that is $w_x \le_\T w_x h$, for all $x \in X$ and $h \in H$. The
latter holds if and only if $w_x$ is the minimum element of $x$ in
$\Abs(W)$ for every coset $x \in X$ and the proof follows.
\qed

\bigskip
A subgroup of $W$ generated by reflections is called a
\textit{reflection subgroup}. The absolute length function on such
a subgroup $K$ is defined with respect to the set of reflections
$\T \cap K$. When $W$ is finite, this function coincides with the
restriction of $\ell_\T: W \to \NN$ on $K$ (this follows from part
(f) of Fact \ref{absolute-length-properties}). As a result, the
corresponding absolute order on $K$ coincides with the induced
order from $\Abs(W)$ on $K$.

\begin{proposition} \label{prop:transitivity}
Assume that $W$ is finite. If $K$ is a modular reflection subgroup
of $W$ and $H$ is a modular subgroup of $K$, then $H$ is a modular
subgroup of $W$.
\end{proposition}

\noindent
\emph{Proof.}
Let $x$ be any left coset of $H$ in $W$. Clearly, $x$ is contained
in a left coset $y$ of $K$ in $W$. Since $K$ is modular in $W$, the
coset $y$ has a minimum element $w_\circ$ in $\Abs(W)$. We leave it
to the reader to check that the map $f: K \mapsto w_\circ K = y$,
defined by $f(w) = w_\circ w$ for $w \in K$, is a poset isomorphism,
where $K$ and $y$ are considered as induced subposets of $\Abs(W)$.
Thus $x$ is isomorphic to its preimage $f^{-1} (x)$ in $K$ under $f$,
which is a left coset of $H$ in $K$. Since $H$ is modular in $K$,
this preimage has a minimum element in $\Abs(K)$, therefore in
$\Abs(W)$, and hence so does $x$. It follows that $H$ is modular
in $W$.
\qed

\begin{remark} \label{remark:transitivity} \rm
The absolute length function on $K$ with respect to $\T \cap K$
coincides with the restriction of $\ell_\T: W \to \NN$ on $K$ even
if $W$ is infinite, provided $K$ is a parabolic reflection
subgroup of $W$ (meaning that $K$ is conjugate to a subgroup
generated by simple reflections) \cite[Corollary 1.4]{Dyer01}.
Thus, the transitivity property of modularity in
Proposition~\ref{prop:transitivity} holds in this situation as
well.
\end{remark}

\begin{proposition} \label{prop:modular-ideal}
Assume $H$ is modular in $W$ and let $\sigma(x)$ be the minimum element
of $x \in X$ in $\Abs(W)$. Then the map $\sigma: X \mapsto W$ induces
a poset isomorphism from $\Abs(X)$ onto an order ideal of $\Abs(W)$.
\end{proposition}

\noindent
\emph{Proof.}
We need to show that (i) $x \le_\T y \Leftrightarrow \sigma (x)
\le_\T \sigma (y)$ for all $x, y \in X$ and that (ii) $\sigma (X)$
is an order ideal of $\Abs(W)$. For $x, y \in X$ we have
  \begin{eqnarray*}
  x\le_\T y & \Leftrightarrow & y=wx \ \, {\rm for \ some} \ \, w \in W
  \ \, {\rm with} \ \, \ell_\T (w) = \ell_\T (y) - \ell_\T (x) \\
  & \Leftrightarrow & w \sigma (x) \in \sigma (y) H \ \, {\rm for \ some}
  \ \, w \in W \ \, {\rm with} \ \, \ell_\T (w) = \ell_\T (\sigma (y)) -
  \ell_\T (\sigma (x)) \\
  & \Leftrightarrow & w \sigma (x) = \sigma (y) \ \, {\rm for \ some}
  \ \, w \in W \ \, {\rm with} \ \, \ell_\T (w) = \ell_\T (\sigma (y)) -
  \ell_\T (\sigma (x)) \\
  & \Leftrightarrow & \sigma(x) \le_\T \sigma(y),
  \end{eqnarray*}

\noindent where the third equivalence is because $\sigma(y)$ is
the unique element of minimum absolute length in its coset and
$\ell_\T (w\sigma(x)) \le \ell_\T (w) + \ell_T (\sigma(x)) =
\ell_\T(\sigma(y))$. This proves (i).

For (ii), given elements $u, w \in W$ with $u \le_\T w$ and $w
\in \sigma(X)$, we need to show that $u \in \sigma(X)$. We set
$v = u^{-1} w$, so that $uv = w$ and $\ell_\T (w) = \ell_\T (u)
+ \ell_\T (v)$. Since $w$ is the minimum element of $w H$ in
$\Abs(W)$, we have $\ell_\T (wh) = \ell_\T (w) + \ell_\T (h)$
for every $h \in H$. Thus, for $h \in H$ we have
  \begin{eqnarray*}
  \ell_\T (uvh) & = & \ell_\T (wh) \ = \ \ell_\T (w) + \ell_\T (h)
  \ = \ \ell_\T (u) + \ell_\T (v) + \ell_\T (h) \\
  \ell_\T (uvh) & = & \ell_\T (uh \cdot h^{-1}vh) \ \le \ \ell_\T
  (uh) + \ell_\T (h^{-1}vh) \ = \ \ell_\T (uh) + \ell_\T (v).
  \end{eqnarray*}
We conclude that $\ell_\T (uh) \ge \ell_\T (u) + \ell_\T (h)$,
hence that $\ell_\T (uh) = \ell_\T (u) + \ell_\T (h)$, for every
$h \in H$. This means that $u$ is the minimum element of $u H$
in $\Abs(W)$, so that $u \in \sigma(X)$, and the proof follows.
\qed

\begin{remark} \label{rem:induced} \rm
Part (i) of the proof of Proposition~\ref{prop:modular-ideal}
shows that $\Abs(X)$ is isomorphic to an induced subposet of
$\Abs(W)$ (moreover, covering relations are preserved). For
that we only needed that each left coset of $H$ in $W$ has a
unique element of minimum absolute length.
\qed
\end{remark}

Next we give a characterization of modularity (which explains
our choice of terminology) for the class of parabolic reflection
subgroups of $W$.

First we need to recall some background and notation on finite
Coxeter groups. Such a group $W$ acts faithfully on a
finite-dimensional Euclidean space $V$ by its standard geometric
representation \cite[\S V.4]{Bou} \cite[\S V.3]{Humphreys}. This
representation realizes $W$ as a group of orthogonal transformations
on $V$ generated by reflections. Let $\Phi$ be a corresponding
root system. For $\alpha \in \Phi$, we denote by $\hH_\alpha$
the linear hyperplane in $V$ which is orthogonal to $\alpha$ and
by $t_\alpha$ the orthogonal reflection in $\hH_\alpha$, so that
$\T = \{ t_\alpha: \alpha \in \Phi \}$. We denote by $\lL_\aA$
the intersection lattice \cite[\S 2.1]{OT} \cite[\S 1.2]{Stanley}
of the Coxeter arrangement $\aA = \{ \hH_\alpha: \alpha \in \Phi
\}$ and by $\lL_W$ the geometric lattice of all linear subspaces
of $V$ (flats) spanned by subsets of $\Phi$, partially ordered by
inclusion. Thus $\lL_\aA$ and $\lL_W$ are isomorphic as posets
and the map which sends an element of $\lL_\aA$ to its orthogonal
complement in $V$ is a poset isomorphism from $\lL_\aA$ onto
$\lL_W$.

Given a reflection subgroup $H$ of $W$, we will denote by $V_H$
the linear span of all roots $\alpha \in \Phi$ for which $t_\alpha
\in H$, so that $V_H \in \lL_W$. Then $H$ is parabolic if and only
if $t_\alpha \in H$ for every $\alpha \in \Phi \cap V_H$ (see, for
instance, \cite{BarceloIhrig99}). Finally, we recall that an
element $Z$ of a geometric lattice $\lL$ is called \emph{modular}
\cite[Definition 2.25]{OT} \cite{Sta71} \cite[Definition
4.12]{Stanley} if we have
  $$ \rk (Y) + \rk(Z) \ = \ \rk (Y \wedge Z) + \rk (Y \vee Z) $$
for every $Y \in \lL$, where $\rk: \lL \mapsto \NN$ denotes the
rank function of $\lL$ and $Y \wedge Z$ (respectively, $Y \vee Z$)
stands for the greatest lower bound (respectively, least upper
bound) of $Y$ and $Z$ in $\lL$.

\begin{theorem} \label{thm:modular}
Assume that $W$ is finite and that $H$ is a parabolic reflection
subgroup of $W$. Then $H$ is a modular subgroup of $W$ if and only
if $V_H$ is a modular element of the geometric lattice $\lL_W$.
\end{theorem}

We will give two proofs of Theorem \ref{thm:modular}. We first
need to establish two crucial lemmas. We recall \cite[\S
2.4]{Armstrong} that to any $w \in W$ are associated the spaces
$\Fix(w) \in \lL_\aA$ and $\Mov(w) \in \lL_W$, where $\Fix(w)$ is
the set of points in $V$ which are fixed by the action of $w$ and
$\Mov(w)$ is the orthogonal complement of $\Fix (w)$ in $V$. For
instance, for every $\alpha \in \Phi$ the space $\Mov(t_\alpha)$
is the one-dimensional subspace of $V$ spanned by $\alpha$. The
maps $\Fix: W \mapsto \lL_\aA$ and $\Mov: W \mapsto \lL_W$ are
surjective and we have $\dim \Mov (w) = \ell_\T (w)$ for every
$w \in W$. Moreover (see the proof of \cite[Theorem
2.4.7]{Armstrong}), if $w = t_{\alpha_1} t_{\alpha_2} \cdots
t_{\alpha_k}$ is a reduced $\T$-word for $w$, then $\{ \alpha_1,
\alpha_2,\dots,\alpha_k \}$ is an $\RR$-basis of $\Mov (w)$. In
particular, $u \le_\T v \Rightarrow \Mov (u) \subseteq \Mov (v)$
for $u, v \in W$.

\begin{lemma} \label{lem:modular}
Assume that $W$ is finite and that $H$ is a reflection subgroup
of $W$ and let $w_\circ \in W$. Then $w_\circ$ is the minimum
of $w_\circ H$ in $\Abs(W)$ if and only if $\Mov (w_\circ) \cap
V_H = \{0\}$.
\end{lemma}

\noindent
\emph{Proof.}
Let $w_\circ = t_{\alpha_1} t_{\alpha_2} \cdots t_{\alpha_k}$ be
a reduced $\T$-word for $w_\circ$. Thus $\ell_\T (w_\circ) = k$
and $\{ \alpha_1, \alpha_2,\dots,\alpha_k \}$ is an $\RR$-basis
of $\Mov (w_\circ)$.

Suppose first that $\Mov (w_\circ) \cap V_H = \{0\}$. We need to
show that $w_\circ \le_\T w_\circ h$ for every $h \in H$. Let $h
= t_{\beta_1} t_{\beta_2} \cdots t_{\beta_\ell}$ be a reduced
$\T$-word for $h \in H$. Then $\ell_\T (h) = \ell$ and $\{ \beta_1,
\beta_2,\dots,\beta_\ell \}$ is an $\RR$-basis of $\Mov (h)$.
Since $h$ is a product of reflections in $H$, its fixed space
contains the orthogonal complement of $V_H$ and hence $\Mov (h)
\subseteq V_H$. We conclude that $\{ \beta_1,
\beta_2,\dots,\beta_\ell \}$ is a linearly independent subset of
$V_H$. Our hypothesis implies that $\{ \alpha_1,\dots,\alpha_k,
\beta_1,\dots,\beta_\ell \}$ is a linearly independent subset of
$V$. We may infer from Carter's Lemma \cite[Lemma 2.4.5]{Armstrong}
that $t_{\alpha_1} \cdots t_{\alpha_k} t_{\beta_1} \cdots
t_{\beta_\ell}$ is a reduced $\T$-word for $w_\circ h$. Therefore
$\ell_\T (w_\circ h) = \ell_\T (w_\circ) + \ell_\T (h)$, which
means that $w_\circ \le_\T w_\circ h$.

Conversely, suppose that $w_\circ$ is the minimum of $w_\circ H$
in $\Abs(W)$. We choose an $\RR$-basis $\{ \beta_1,
\beta_2,\dots,\beta_\ell \}$ of $V_H$ consisting of roots $\beta_i$
with $t_{\beta_i} \in H$ and set $h = t_{\beta_1} t_{\beta_2}
\cdots t_{\beta_\ell} \in H$. By assumption, we have $\ell_\T
(w_\circ h) = \ell_\T (w_\circ) + \ell_\T (h)$. This equation and
Carter's Lemma imply that $\{ \alpha_1,\dots,\alpha_k,
\beta_1,\dots,\beta_\ell \}$ is linearly independent or,
equivalently, that $\Mov (w_\circ) \cap V_H = \{0\}$.
\qed

\begin{lemma} \label{lem:minimal}
Assume that $W$ is finite and that $H$ is a parabolic reflection
subgroup of $W$ and let $w \in W$. Then $w$ is a minimal element
of $wH$ in $\Abs(W)$ if and only if $\Mov (w) \wedge V_H = \{0\}$
holds in $\lL_W$.
\end{lemma}

\noindent
\emph{Proof.}
We recall that every element of $\lL_W$ is of the form $\Mov (u)$
for some $u \in W$ and that $\Mov (u)$ is nonzero if and only if it
contains $\Mov (t)$ for some $t \in \T$. Moreover, we have $\Mov (t)
\subseteq \Mov (u) \Leftrightarrow t \le_\T u$ \cite[Theorem
2.4.7]{Armstrong} for $t \in \T$ and since $H$ is parabolic, we have
$t \in H$ for every reflection $t \in \T$ for which $\Mov (t)
\subseteq V_H$. From these facts we conclude that $\Mov (w) \wedge
V_H \ne \{0\}$ holds in $\lL_W$ if and only if there exists $t \in H
\cap \T$ such that $t \le_\T w$. The latter holds if and only if $wt
<_\T w$ for some $t \in H \cap \T$ or, equivalently, if and only if
$w$ is not a minimal element of $wH$ in $\Abs(W)$.
\qed

\bigskip
\noindent
\textit{First proof of Theorem \ref{thm:modular}.} We will use the
following characterization of modularity in $\lL_W$: An element $Z \in
\lL_W$ is modular if and only if $Y \cap Z \in \lL_W$ for every $Y
\in \lL_W$. This statement follows directly from \cite[Lemma 2.24]{OT},
which implies that an element $Z \in \lL_\aA$ is modular if and only
if $Y + Z \in \lL_\aA$ for every $Y \in \lL_\aA$.

We first assume that $H$ is modular in $W$ and consider any element
$Y \in \lL_W$. We need to show that $Y \cap V_H \in \lL_W$. Since
$Y \in \lL_W$, we have $Y = \Mov (w)$ for some $w \in W$. By our
assumption, the coset $wH$ has a minimum element, say $w_\circ$, in
$\Abs(W)$. We claim that $Y \cap V_H = \Mov (w_\circ^{-1} w)$. Since
$\Mov (w_\circ^{-1} w) \in \lL_W$, it suffices to prove the claim.
Indeed, since $w_\circ \le_\T w$, we also have $w_\circ^{-1} w \le_\T
w$ and hence $\Mov (w_\circ^{-1} w) \subseteq \Mov (w) = Y$. Similarly,
since $w \in w_\circ H$, we have $w_\circ^{-1} w \in H$ and hence
$\Mov (w_\circ^{-1} w) \subseteq V_H$, so we may conclude that $\Mov
(w_\circ^{-1} w) \subseteq Y \cap V_H$. For the reverse inclusion,
we recall \cite[p.~25]{Armstrong} that
  $$ Y \ = \ \Mov (w) \ = \ \Mov (w_\circ) \, \oplus \, \Mov
     (w_\circ^{-1} w). $$
By our choice of $w_\circ$ and Lemma \ref{lem:modular} we have $\Mov
(w_\circ) \cap V_H = \{0\}$. As we already know that $Y \cap V_H
\supseteq  \Mov (w_\circ^{-1} w)$, it follows that $Y \cap V_H = \Mov
(w_\circ^{-1} w)$.

Suppose now that $V_H$ is a modular element of $\lL_W$ and consider
any left coset $x$ of $H$ in $W$. We need to show that $x$ has a
minimum in $\Abs(W)$. Let $w_\circ$ be any minimal element of $x$ in
$\Abs(W)$. Since $\Mov (w_\circ) \cap V_H \in \lL_W$, by modularity
of $V_H$, the greatest lower bound $\Mov (w_\circ) \wedge V_H$ of
$\Mov (w_\circ)$ and $V_H$ in $\lL_W$ must be equal to $\Mov (w_\circ)
\cap V_H$. This statement and Lemmas \ref{lem:modular} and
\ref{lem:minimal} imply that $w_\circ$ is the minimum element of $x$
in $\Abs(W)$ and the proof follows.
\qed

\begin{remark} \label{rem:needparabolic} \rm
The assumption in Theorem \ref{thm:modular} that the reflection
subgroup $H$ is parabolic was not used in the proof of the \emph{only
if} direction of the theorem. However, it is essential for the other
direction. Indeed, let $W$ be the dihedral group of symmetries of a
square $Q$ and let $H$ be the subgroup of order 4 generated by the
reflections on the lines through the center of $Q$ which are
parallel to the sides. The unique left coset of $H$ in $W$, other
than $H$, has no minimum element in $\Abs(W)$ and hence $H$ is not
modular in $W$. On the other hand, $V_H = V$ is trivially a modular
element of the lattice $\lL_W$.
\qed
\end{remark}

\medskip
For the second proof of Theorem \ref{thm:modular} we recall the
following definition. Let $\lL$ be a geometric lattice of rank $d$,
with rank function $\rk: \lL \mapsto \NN$. The \emph{characteristic
polynomial} of $\lL$ is defined by the formula
  \begin{equation} \label{eq:chardef}
    \chi_\lL (q) \ := \ \sum\limits_{Y \in \lL} \, \mu_\lL
    (\hat{0}, Y) \, q^{d - \rk(Y)},
  \end{equation}
where $\mu_\lL$ stands for the M\"obius function \cite[\S 3.7]{StaEC1}
of $\lL$ and $\hat{0}$ is the minimum element of $\lL$. We now let
$\lL = \lL_W$ and recall that $\rk(Y) = \dim (Y)$ and (see, for instance,
\cite[Lemma 4.7]{OS80})
  \begin{equation} \label{eq:mobiusW}
    (-1)^{\rk(Y)} \mu_\lL (\hat{0}, Y) \ = \ \# \{w \in W: \Mov (w) = Y\}
  \end{equation}
for $Y \in \lL$ and that $\dim \Mov (w) = \ell_\T (w)$ for $w \in W$. As
a result, the characteristic polynomial of $\lL_W$ is related to the rank
generating polynomial of $\Abs(W)$ by the well known equality
  \begin{equation}\label{eq:charabs}
    W_\T (q) \ = \ (-q)^d \, \chi_\lL (-1/q).
  \end{equation}

\bigskip
\noindent
\textit{Second proof of Theorem \ref{thm:modular}.} Let us write $\lL =
\lL_W$, as before, and set $Z = V_H \in \lL$. By the Modular Factorization
Theorem for geometric lattices \cite{Sta71} \cite[Theorem 4.13]{Stanley}
and its converse (see \cite[Section 8]{Ku93}) we have that $Z$ is a
modular element of $\lL$ if and only if
 \begin{equation}\label{eq:MFT}
    \chi_\lL (q) \ = \ \chi_{[\hat{0}, Z]} (q) \ \left(
    \sum\limits_{Y \in \lL: \ Y \wedge Z = \hat{0}} \, \mu_\lL
    (\hat{0}, Y) \, q^{d - \rk(Y) - \rk(Z)} \right),
  \end{equation}
where $[\hat{0}, Z]$ denotes a closed interval in $\lL$ and $\hat{0} =
\{0\}$ is the minimum element of $\lL$. Replacing $q$ by $-1/q$ and
taking (\ref{eq:mobiusW}) and (\ref{eq:charabs}) into account, we see
that (\ref{eq:MFT}) can be rewritten as
  \begin{equation}\label{eq:MFT3}
    W_\T (q) \ = \ H_\T (q) \ \left( \sum\limits_{\Mov (w) \wedge Z
    = \hat{0} } \ q^{\ell_\T (w)} \right).
  \end{equation}
We recall that every finite partially ordered set has at least one
minimal element. Assume first that $Z$ is a modular element of $\lL$.
Setting $q=1$ in (\ref{eq:MFT3}) and using Lemma \ref{lem:minimal} we
conclude that every left coset of $H$ in $W$ has exactly one minimal
(and hence a minimum) element in $\Abs(W)$. By definition, this means
that $H$ is a modular subgroup of $W$. Conversely, suppose that $H$ is
a modular subgroup of $W$. Then, by Lemma \ref{lem:minimal}, the sum
in the right-hand side of (\ref{eq:MFT3}) is equal to $X_\T (q)$ and
hence (\ref{eq:MFT3}) holds by Proposition \ref{prop:fact}. Thus $Z$
is a modular element of $\lL$ and the proof follows.
\qed

\begin{proposition} \label{prop:refmodular}
Assume that $W$ is finite. Then every modular reflection subgroup
of $W$ is a parabolic reflection subgroup.
\end{proposition}

\noindent
\emph{Proof.}
Let $H$ be a modular reflection subgroup of $W$ and let $K$ be the
unique parabolic reflection subgroup of $W$ with $V_K = V_H$. Thus
$K$ is generated by all reflections $t \in \T$ with $\Mov (t)
\subseteq V_H$ and contains $H$ as a reflection subgroup. We need
to show that $H = K$. Since $H$ is modular in $W$, it is also modular
in $K$. Thus, without loss of generality we may assume that $K = W$,
so that $\rank(H) = \rank(W)$. We note that $H_\T (q)$ and $W_\T(q)$
are both polynomials of degree $\rank (W)$. Therefore, Proposition
\ref{prop:fact} implies that $X_\T (q)$ is a constant. Since this
can only happen if $X$ is a singleton, we conclude that $H = W$ and
the proof follows.
\qed

\begin{question} \label{que:refmodular}
Does there exist a modular subgroup of a Coxeter group which is not
a reflection subgroup?
\end{question}

We recall that a poset $P$ is said to be graded of rank $d$ if
every maximal chain in $P$ has exactly $d+1$ elements. The
following proposition generalizes the fact that $\Abs(W)$ is
graded with rank equal to $\rank(W)$.

\begin{proposition} \label{prop:graded}
The order $\Abs(X)$ is graded of rank $\rank (W) - \rank(H)$ for
every finite Coxeter group $W$ and every modular reflection subgroup
$H$ of $W$.
\end{proposition}

\noindent
\emph{Proof.}
Since $\Abs(X)$ has a minimum element and is locally graded with rank
function given by absolute length (Proposition \ref{prop:modular-ideal}),
it suffices to show that for every element $x \in X$ there exists $y
\in X$ of absolute length $\rank (W) - \rank(H)$ such that $x \le_\T y$.

Consider any $x \in X$ and let $u_\circ$ be the minimum element of
$x$ in $\Abs(W)$. Thus we have $\Mov (u_\circ) \cap V_H = \{0\}$
by Lemma \ref{lem:modular} and $\ell_\T (x) = \ell_\T (u_\circ) =
\dim \Mov (u_\circ)$. Let $u_\circ = t_{\alpha_k} \cdots
t_{\alpha_2} t_{\alpha_1}$ be a reduced $\T$-word for $u_\circ$,
so that $\ell_\T (x) = k$. We extend $\{ \alpha_1,
\alpha_2,\dots,\alpha_k \}$ to a maximal linearly independent set
of roots $\{ \alpha_1, \alpha_2,\dots,\alpha_r\}$ whose linear
span intersects $V_H$ trivially and set $w_\circ = t_{\alpha_r}
\cdots t_{\alpha_2} t_{\alpha_1}$ and $y = w_\circ H \in X$.
Clearly, we have $r = \dim (V) - \dim (V_H) = \rank (W) - \rank
(H)$. Since $\Mov (w_\circ)$ is the linear span of $\alpha_1,
\alpha_2,\dots,\alpha_r$, we have $\Mov (w_\circ) \cap V_H =
\{0\}$ by construction. Lemma \ref{lem:modular} implies that
$w_\circ$ is the minimum element of $y$ in $\Abs(W)$ and hence
that $\ell_\T (y) = \ell_\T (w_\circ) = r$. Finally, setting $v =
w_\circ u_\circ^{-1} = t_{\alpha_r} \cdots t_{\alpha_{k+1}}$ we
have $w_\circ = v u_\circ$ and hence $y = vx$. By Carter's Lemma
\cite[Lemma 2.4.5]{Armstrong} we also have $\ell_\T (v) = r-k =
\ell_\T (y) - \ell_\T (x)$. Definition \ref{def:main1} implies
that $x \le_\T y$ and the proof follows.
\qed

\begin{question} \label{que:graded}
Does there exist a subgroup $H$ of a Coxeter group $W$ for which
$\Abs(X)$ is not graded?
\end{question}

A reflection subgroup $H$ of $W$ is said to be of \emph{almost
maximal rank} if $\rank(H) = \rank(W) - 1$. Modular parabolic
reflection subgroups of this kind can be characterized as follows.

\begin{proposition} \label{prop:maxmodular}
Assume that $W$ is finite and that $H$ is a parabolic reflection
subgroup of $W$, other than $W$. The following are equivalent:
\begin{itemize}
\itemsep=0pt
\item[{\rm (i)}]
$H$ is a modular subgroup of $W$ of almost maximal rank.
\item[{\rm (ii)}]
Every left coset of $H$, other than $H$, contains a reflection.
\item[{\rm (iii)}]
Every left coset of $H$, other than $H$, contains a unique
reflection.
\end{itemize}
\end{proposition}

\noindent
\emph{Proof.}
Suppose that (i) holds. We then have $W_\T(q) = H_\T(q) \, X_\T (q)$
by Proposition \ref{prop:fact}. Since the degrees of $W_\T(q)$ and
$H_\T(q)$ are equal to the Coxeter ranks of $W$ and $H$, respectively,
it follows that the degree of $X_\T (q)$ is equal to one. This means
that every left coset $x \in X$ of $H$, other than $H$, contains an
element of absolute length one, so that (ii) is satisfied. We have
shown that (i) implies (ii).

Suppose that (ii) holds and let $x \in X$ be a left coset of $H$ in
$W$, other than $H$. Choose a reflection $t \in x$. Since $H$ is
parabolic and does not contain $t$, we have $\Mov (t) \cap V_H =
\{0\}$. Lemma \ref{lem:modular} implies that $t$ is the minimum element
of $x$ in $\Abs(W)$. In particular, $x$ contains a unique reflection.
We conclude that (ii) implies both (i) and (iii). The implication
${\rm (iii)} \Rightarrow {\rm (ii)}$ is trivial.
\qed

\begin{question} \label{que:almostmax}
Does there exist a non-parabolic (necessarily non-modular) reflection
subgroup $H$ of a finite Coxeter group  such that every left coset,
other than $H$, contains a unique reflection?
\end{question}

\begin{example} \label{ex:SkBk} \rm
For $k \le n$ and under the natural embedding, the symmetric and
hyperoctahedral groups $S_k$ and $B_k$ are modular subgroups of
$S_n$ and $B_n$, respectively. This follows from Theorem
\ref{thm:modular} and known facts on the modular elements of the
geometric lattice $\lL_W$ in these cases; see, for instance,
\cite[Theorem 2.2]{BarceloIhrig98}. Alternatively, one can check
directly that for $1 \le i \le n-1$, the transpositions $(i \ n)$
are representatives of the left cosets of $S_{n-1}$ in $S_n$,
other than $S_{n-1}$. Proposition \ref{prop:maxmodular} implies
that $S_{n-1}$ is modular in $S_n$. The transitivity property of
Proposition \ref{prop:transitivity} implies that $S_k$ is modular
in $S_n$ for each $k \le n$. A similar argument works for the
hyperoctahedral groups. \qed
\end{example}

We end this section with two more open questions.

\begin{question} \label{que:infinite}
Do infinite modular subgroups exist?
\end{question}

\begin{question} \label{que:max}
For which subgroups $H$ of $W$ does $\Abs(X)$ have a maximum
element?
\end{question}

\section{Quasi-modular subgroups}
\label{sec:quasi}

This section introduces a condition on a subgroup of a Coxeter
group, termed quasi-modularity, which is broader than modularity
and guarantees an affirmative answer to Question \ref{que:fact}.
Examples of quasi-modular subgroups which are not modular are
discussed. Throughout this section, the set of reflections of a
Coxeter group $H$ will be denoted by $\T(H)$.

\subsection{Quasi-modularity}
\label{subsec:quasi}

The main definition of this section is as follows.

\begin{defn} \rm
A subgroup $H$ of a finite Coxeter group $W$ is {\em quasi-modular}
if $H$ is isomorphic to a Coxeter group and
\begin{equation} \label{eq:quasi}
W_\T (q) \ = \ H_{\T(H)} (q) \cdot X_\T (q),
\end{equation}
where $\T = \T(W)$ and $\T(H)$ is the subset of $H$ which corresponds
to the set of reflections of this Coxeter group.
\end{defn}

Proposition~\ref{prop:fact}
implies that for reflection subgroups of $W$, quasi-modularity is
equivalent to modularity. However, this is not the case for general
subgroups as $\T(H)$ may not be equal to $H \cap \T(W)$.

\begin{example} \rm
We list two families of examples of quasi-modular subgroups which
are not modular.

(a) Let $W$ be the Weyl group of type $D_n$, considered as a group
of signed permutations of $\{1, 2,\dots,n\}$ with an even number of
sign changes. Let $H$ be the subgroup consisting of all $w \in
W$ satisfying $w (n) \in \{n, -n\}$. Then $H$ is isomorphic to the
hyperoctahedral group $B_{n-1}$ and the identity element $e \in W$
together with the reflections $( \!(i \ \, n) \! )$ for $1 \le i
\le n-1$ (where the notation is as in Example \ref{ex:modular1} (d))
form a complete list of coset representatives of $H$ in
$W$. As a result, we have $X_\T (q) = 1 + (n-1)q$, where $X = W/H$
and $\T = \T(W)$. Using this fact and (\ref{eq:rankgenW}), it can
be easily verified that (\ref{eq:quasi}) holds in this situation
and hence that $H$ is a quasi-modular subgroup of $W$. On the
other hand, it is also easy to verify that $H_\T (q)$ has degree
$n$, as does $W_\T (q)$. Thus $H$ is not a modular subgroup of $W$
by Proposition \ref{prop:fact}.

(b)  Consider the symmetric group $S_{2n}$ as the group of all
permutations of the set $\Omega_n := \{1, -1, 2, -2,\dots,n, -n\}$
and the natural embedding of the hyperoctahedral group $B_n$ in $S_{2n}$,
mapping the Coxeter generators of $B_n$ to the transposition $(n \,
-n)$ and the products $(i \ i+1) (-i \ -i-1)$ for $1 \le i \le n-1$.
Clearly, this embedded copy of $B_n$ is not a reflection subgroup of
$S_{2n}$. Several combinatorial interpretations to the poset
$\Abs(S_{2n}/B_n)$ will be given in Section \ref{subsection:PM},
where the following statement will also be proved.

\begin{theorem} \label{pm-quasi-main}
The group $B_n$ is a non-modular, quasi-modular subgroup of $S_{2n}$
for every $n \ge 2$.
\end{theorem}

\end{example}

\subsection{Balanced complex reflections}

Before proving Theorem~\ref{pm-quasi-main} we introduce an absolute
order on balanced complex reflections. Recall that the wreath product
of the cyclic group $\ZZ_r$ by the symmetric group $S_n$ is defined
as
\[ G(r,n) \ = \ \ZZ_r \wr S_n \ := \
\{ [(c_1,\ldots, c_n);\,\pi]: \ c_i \in \ZZ_r, \, \pi \in S_n\} \]
with group operation
\[ [(c_1,\ldots, c_n); \, \pi] \cdot [(c'_1,\dots, c'_n);\,\pi'] \
:= \ [(c_1+c'_{\pi^{-1}(1)},\dots,c_n+c'_{\pi^{-1}(n)});\,\pi\pi']. \]

We think of the elements of $\ZZ_r$ as colors and denote by $\psi:
G(r,n) \to S_n$ the canonical map, defined by $\psi( [\bar c; \,
\pi]) := \pi$. Via this map, the elements of $G(r, n)$ inherit a
cycle structure from those of $S_n$.

\begin{defn} \rm
A cycle of an element of $G(r,n)$ is {\em balanced} if the sum of
the colors of its elements is zero modulo $r$. An element $w \in
G(r,n)$ is {\em balanced} if all cycles of $w$ are balanced. We
denote by $C(r,n)$ the set of balanced elements of $G(r,n)$.
\end{defn}

For example, there are three balanced elements in $G(2,2) \cong B_2$,
namely the identity and the reflections $[(0,0);\,(1 \ 2)] = (1 \
2)(-1 \, -2)$ and $[(1,1);\,(1 \ 2)] = (1 \, -2) (-1 \ 2)$. Note
that $C(2, 2)$ is not a subgroup of $G(2, 2)$.

\begin{remark} \rm
To motivate the notion of balanced, we note that balanced cycles
generalize the notion of {\em paired cycles}, introduced by Brady
and Watt~\cite{BW} in the study of the absolute order of types $B$
and $D$ and further studied in~\cite{Ka}. Moreover, conjugacy
classes in $G(r, n)$ are parametrized by cycle type and sum of
colors (modulo $r$) in each cycle (so that $C(r,n)$ is the union
of conjugacy classes in $G(r, n)$).
%
\end{remark}

The wreath product $G(r,n)$ acts on the vector space $V = \CC^n$
by permuting coordinates and multiplying them by suitable $r$th
roots of unity, in a standard way. The set of pseudoreflections
$\T(r,n)\subseteq G(r,n)$ consists of all elements fixing a
hyperplane (codimension one subspace). The absolute length
function $\ell_{\T(r,n)}: G(r, n) \to \NN$ is defined with
respect to the generating set $\T(r,n)$.

\begin{remark} \label{rem:balanced} \rm
(a) We have $\psi(\T(r, n)) = \T (S_n) \cup \{e\}$, where $e \in
S_n$ stands for the identity element.

(b) The set $\T(r,n) \cap C(r,n)$ of balanced pseudoreflections
in $G(r, n)$ consists of all elements of the form $\tau = [\bar{c};
\, t]$, where $t = (a \ b) \in \T (S_n)$ and $\bar{c}$ assigns
opposite colors to $a$ and $b$ and the zero color to all other
elements of $\{1, 2,\dots,n\}$. As a result, we have $\psi(\T(r,
n) \cap C(r,n)) = \T (S_n)$.

(c) The canonical map $\psi$ has the following crucial property:
given $w \in C(r,n)$ and $t \in \T(S_n)$ such that $t \psi(w)$
is covered by $\psi(w)$ in $\Abs(S_n)$, there is a unique
(necessarily balanced) pseudoreflection $\tau \in \psi^{-1}(t)$
such that $\tau w \in  C(r, n)$.
\end{remark}

\begin{defn} \rm
The absolute order on $C(r,n)$, denoted $\Abs(C(r,n))$, is the
reflexive and transitive closure of the relation consisting of
the pairs $(u, v)$ of elements of $C(r,n)$ for which $v = \tau u$
for some $\tau \in \T(r,n) \cap C(r,n)$ and $\ell_{\T(r,n)}(u) <
\ell_{\T(r,n)}(v)$.
\end{defn}

The partial order $\Abs(C(r,n))$ is the subposet induced on $C(r,n)$
from Shi's absolute order on $G(r,n)$; see~\cite{Shi1, Shi2}. We
will focus on this subposet since it will be useful (in the special
case $r=2$) in our proof of Theorem~\ref{pm-quasi-main}.

\begin{proposition}\label{obs-canonical}
\begin{itemize}
\itemsep=0pt
\item[{\rm (a)}] The canonical map $\psi: G(r,n) \to S_n$ induces
a rank preserving poset epimorphism from the order $\Abs (C(r,n))$
onto $\Abs(S_n)$.

\item[{\rm (b)}] Every maximal interval in $\Abs(C(r, n))$ is mapped
isomorphically by $\psi$ onto a maximal interval in $\Abs(S_n)$.
\end{itemize}
\end{proposition}

\noindent
\emph{Proof.}
(a) The map $\psi$ is a group epimorphism, by its definition, and
$\psi(\T(r,n)) = \T(S_n) \cup \{e\}$ by Remark~\ref{rem:balanced}
(a). Hence we have $\ell_{\T(r,n)} (w) \ge \ell_\T (\psi(w))$ for
every $w \in G(r,n)$. For $w \in C(r,n)$ the reverse inequality
$\ell_{\T(r,n)} (w) \le \ell_\T (\psi(w))$ follows from
Remark~\ref{rem:balanced} (c). Thus we have $\ell_{\T(r,n)} (w) =
\ell_\T (\psi(w))$ for every $w \in C(r,n)$. Furthermore, this fact
and parts (b) and (c) of Remark~\ref{rem:balanced} imply that for
$u, v \in C(r,n)$, we have $v = \tau u$ for some $\tau \in \T(r,n)
\cap C(r,n)$ and $\ell_{\T(r,n)} (u) < \ell_{\T(r,n)}(v)$ if and
only if $\psi (v) = t\psi(u)$ for some $t \in \T(S_n)$ and $\ell_\T
(\psi(u)) < \ell_\T (\psi(v))$. In other words, $u$ is covered by
$v$ in $\Abs(C(r,n))$ if and only if $\psi(u)$ is covered by $\psi
(w)$ in $\Abs(S_n)$.

(b) We first check that $\psi$ maps maximal elements of $\Abs(C(r,
n))$ to maximal elements of $\Abs(S_n)$. Indeed, since $\psi$ is rank
preserving, the rank of an element $w$ in $\Abs(C(r,n))$ is equal to
$n - k$, where $k$ is the number of cycles. Thus, if $\psi(w)$ is
not maximal in $\Abs(S_n)$, then $\psi(w)$ has at least two cycles and
one can check that there exists $\tau \in \T(r,n) \cap C(r, n)$ such
that $\tau w$ has fewer cycles than $w$, so $w$ is not maximal
either. We next observe that by Remark~\ref{rem:balanced} (c), for
every $w \in C(r,n)$ the map $\psi$ induces a bijection between
elements covered by $w$ in $\Abs (C(r,n))$ and those covered by $\psi
(w)$ in $\Abs(S_n)$. By induction on the rank of the top element, it
follows that intervals in $\Abs (C(r,n))$ are mapped isomorphically
by $\psi$ to intervals in $\Abs(S_n)$. In particular, every maximal
interval in $\Abs(C(r,n))$ is mapped isomorphically by $\psi$ onto a
maximal interval in $\Abs(S_n)$.
\qed

\begin{corollary} \label{rank-function-rn}
\[ \sum\limits_{w\in C(r,n)} q^{\ell_{\T(r,n)}(w)} \ = \
\prod\limits_{i=1}^{n-1} \, (1+ riq). \]
\end{corollary}

\noindent
\emph{Proof.}
Let $\psi_0: C(r,n) \to S_n$ be the restriction of $\psi$ to $C(r,n)$.
By Proposition~\ref{obs-canonical} we have $\ell_{\T(r,n)}(w) =
\ell_\T(\pi) = n-k$ for every $w \in C(r,n)$, where $k$ is the number
of cycles of $\pi := \psi_0(w)$. Since all elements in the preimage
$\psi^{-1}_0 (\pi)$ are balanced, we have
$$
|\psi^{-1}_0 (\pi)| \ = \ r^{n-k} \ = \ r^{\ell_T(\pi)}
 $$
and thus

\begin{eqnarray*}
\sum\limits_{w\in C(r,n)} q^{\ell_{\T(r,n)}(w)} &=&
\sum\limits_{\pi\in S_n} |\psi^{-1}_0 (\pi)| \, q^{\ell_{\T}(\pi)} \ = \
\sum\limits_{\pi\in S_n} r^{\ell_T(\pi)} q^{\ell_{\T}(\pi)} \\
& & \\
&=& \sum\limits_{\pi \in S_n} (rq)^{\ell_{\T}(\pi)} \ = \
\prod \limits_{i=1}^{n-1} \, (1 + riq).
\end{eqnarray*}
\qed


\subsection{Perfect matchings}
\label{subsection:PM}

A partition of set $\Omega$ into two-element subsets is called a
{\it perfect matching}. Throughout this section we will denote by
$\mM_n$ the set of perfect matchings of $\Omega_n = \{1, -1, 2,
-2,\dots,n, -n\}$. Consider the simple graph $\Delta_n$, introduced
in~\cite{HHN}, on the set of nodes $\mM_n$ in which two perfect
matchings are adjacent if their symmetric difference is a cycle of
length 4. The diameter and the enumeration of geodesics of this
graph were studied in~\cite{Ben-Ari}; the induced subgraph on
non-crossing perfect matchings was studied earlier in~\cite{HHN}.

\begin{defn} \rm
Fix an arbitrary element $x_0 \in \mM_n$. The {\it absolute order}
on $\mM_n$, denoted $\Abs(\mM_n)$, is the poset $(\mM_n, \preceq)$
defined by letting $x \preceq y$ if $x$ lies in a geodesic path in
$\Delta_n$ with endpoints $x_0$ and $y$, for $x, y \in \mM_n$.
\end{defn}

The symmetric group $S_{2n}$ of permutations of $\Omega_n$ acts
naturally on $\mM_n$ (this action may be identified with the
conjugation action of $S_{2n}$ on the set of fixed point free
involutions on a $2n$-element set). The stabilizer of $x_0 = \{
\{-1, 1\}, \{-2, 2\},\dots,\{-n, n\} \}$ is the natural embedding
of the hyperoctahedral group $B_n$ in $S_{2n}$ and hence we get
the following statement.

\begin{observation} \label{obs-matching-so}
The poset $\Abs(\mM_n)$ is isomorphic to $\Abs(S_{2n}/B_n)$.
\end{observation}

In particular, the isomorphism type of $\Abs(\mM_n)$ is independent
of the choice of $x_0$. Without loss of generality, for the
remainder of this section we will assume that $x_0$ consists
of the sets (arcs) $\{-i, i\}$ for $1 \le i \le n$.

\begin{proposition}\label{matching-iso}
The poset $\Abs(\mM_n)$ is isomorphic to $\Abs(C(2,n))$.
\end{proposition}

\noindent
\emph{Proof.}
The proof generalizes a construction from~\cite{HHN}.

Given a perfect matching $x \in \mM_n$, consider the union $x \cup
x_0$, consisting of the arcs of $x$ and $\{-i, i\}$ for $1 \le i
\le n$. This is a disjoint union of nontrivial cycles and isolated
arcs. We orient the nontrivial cycles in the following way: Given
any such cycle $C$, we let $k$ be the minimum positive integer
such that $\{-k, k\}$ is an arc of $C$ and choose the cyclic
orientation of $C$ in which this edge is directed from $-k$ to
$k$. We associate to $x$ a signed permutation $f(x) = \pi \in B_n$
as follows. For $i \in \Omega_n$, we set $\pi (i) = i$ if $\{-i,
i\} \in x$. Otherwise we set $\pi(i) = -j$ if either $(i, j)$ or
$(-i,-j)$ is a directed edge in the above orientation, and $\pi(i)
= j$ if either $(-i, j)$ or $(i,-j)$ is a directed edge in the
orientation. We will show that $f: \mM_n \to C(2,n)$ is a
well-defined map which is an isomorphism of the corresponding
absolute orders.

%

We first observe that the map $f: \mM_n \to B_n$ is well-defined.
Indeed, this holds since $\{-i, i\} \in x \cup x_0$ for $1 \le i
\le n$ and hence at most one of $i$ and $-i$ can be the initial
vertex of a directed arc in the above orientation. Moreover,
since the number of arcs of any nontrivial cycle of $x \cup x_0$
is even, the number of arcs with vertices of same sign in such a
cycle must also be even. This implies that every nontrivial cycle
of the signed permutation $f(x)$ is balanced and hence we have
a well-defined map $f: \mM_n \to C(2,n)$.

To show that $f: \mM_n \to C(2,n)$ is a bijection, it suffices
to describe the inverse map $g: C(2,n) \to \mM_n$. Given a balanced
signed permutation $\pi \in C(2,n)$, we construct $g(\pi) \in \mM_n$
as follows. First, we include in $g(\pi)$ the arc $\{-i, i\}$ for
each $i \in \Omega_n$ with $\pi (i) = i$. Second, let $(a_1 \ a_2
\, \cdots \, a_k)$ be any nontrivial cycle of $\pi$ and assume that
$a_1$ is the minimum of the absolute values of the element of this
cycle. We then include in $g(\pi)$ the arcs $\{a_1, -a_2\}$, $\{a_2,
-a_3\},\dots,\{a_k, -a_1\}$. We leave it to the reader to verify
that $g$ is the inverse map of $f$.

Finally we prove that $f: \mM_n \to C(2,n)$ induces an
isomorphism of absolute orders. We consider the simple graph
$\Gamma_n$ on the node set $C(2,n)$ in which two permutations
$\pi, \sigma \in C(2,n)$ are adjacent if $\pi^{-1}\sigma \in
\T(2,n)$. Since $f$ maps $x_0$ to the identity element of $C(2,n)$,
it suffices to show that $f$ induces a graph isomorphism from
$\Delta_n$ to $\Gamma_n$. Indeed, two matchings $x_1, x_2 \in
\mM_n$ are adjacent in $\Delta_n$ if and only if there exist four
distinct elements $i, j, k, l \in \Omega_n$ such that $x_1
\setminus \{ \{i,j\}, \{k,l\} \} = x_2 \setminus \{ \{i,k\},
\{j,l\} \}$. Without loss of generality, we may assume that
$(i, j)$ and $(k, l)$ are directed edges in the orientation of
$x_1 \cup x_0$. By considering the eight cases determined by the
signs of $i, j, k, l$, one can verify that this happens if and
only if there exists a reflection $\tau \in \{(j,k), (-j,-k)\}
\subseteq \T(2,n)$ such that $f(x_2) = \tau f(x_1)$ and the
proof follows.
\qed

\begin{corollary}\label{core-PM}
There is a poset epimorphism from $\Abs(\mM_n)$ to $\Abs(S_n)$ which
maps every maximal interval in $\Abs(\mM_n)$ isomorphically onto a
noncrossing partition lattice of type $A_{n-1}$.
\end{corollary}

\noindent
\emph{Proof.}
This follows from Propositions~\ref{matching-iso}
and~\ref{obs-canonical} and the fact that every maximal interval in
$\Abs(S_n)$ is isomorphic to the lattice of noncrossing partitions
of the set $\{1, 2,\dots,n\}$.
\qed

\bigskip
The previous corollary implies~\cite[Corollaries 1.6 and 2.2]{HHN}
and~\cite[Theorem 3.20]{Ben-Ari}.

\begin{corollary}\label{pm-cor5}
For every $n \ge 1$ we have
\[ {(\mM_n)}_\T (q) \ = \ \prod \limits_{i=0}^{n-1} \,
(1 + 2iq). \]
\end{corollary}

\noindent
\emph{Proof.}
This follows from Proposition~\ref{matching-iso} and
Corollary~\ref{rank-function-rn}.
\qed

\bigskip
\noindent
{\it Proof of Theorem~\ref{pm-quasi-main}.} That $B_n$ is a
quasi-modular subgroup of $S_{2n}$ follows from
Observation~\ref{obs-matching-so}, Corollary~\ref{pm-cor5} and
the known formulas for the rank generating functions of
$\Abs(S_{2n})$ and $\Abs(B_n)$.

Suppose that $B_n$ were a modular subgroup of $S_{2n}$ for some
$n \ge 2$. Then, according to Proposition~\ref{prop:fact} and
Observation~\ref{obs-matching-so} we should have ${(S_{2n})}_\T (q)
= {(B_n)}_\T (q) \cdot {(\mM_n)}_\T (q)$, where $\T := \T(S_{2n})$,
and hence ${(B_n)}_\T (q)$ should have degree $n$. This is not
correct, since there exist elements of the natural embedding of $B_n$
in $S_{2n}$ which are cycles in $S_{2n}$ of absolute length $2n-1$.
\qed

\begin{remark} \rm
By Corollary~\ref{pm-cor5}, the $S_{2n}$ conjugation action on fixed
point free involutions of $\Omega_n$ has a quasi-modular stabilizer.
For $1 \le k < n$ such that $n-k$ is even, however, the $S_n$
conjugation action on involutions of $\{1, 2,\dots,n\}$ with $k$
fixed points has a nicely factorized rank generating function even
though its stabilizer is not quasi-modular.
\end{remark}

\begin{question}
For $r \ge 3$, is there a Coxeter group action whose associated
absolute order is isomorphic to $\Abs(C(r,n))$?
\end{question}

The cardinality of $C(r, n)$ is equal to the product
$\prod_{i=1}^{n-1} (1+ri)$ (by Corollary~\ref{rank-function-rn})
and hence to the number of $(r+1)$-ary increasing trees of order
$n$; see, for instance,~\cite{Sloane}.

\begin{question}
For $r \ge 1$, is there a (natural) Coxeter group action on these
trees whose associated absolute order is isomorphic to $\Abs(C(r,n))$?
\end{question}

\section{An application to alternating subgroups}
\label{sec:alt}

Throughout this section $(W,S)$ will be a Coxeter system with set
of reflections $\T = \{w s w^{-1}: w \in W, \, s \in S\}$. The {\it
alternating subgroup} $W^+$ is defined as the kernel of the {\it sign
character} on $W$, which maps every element of $S$ to $-1$. We will
show that a natural absolute order on $W^+$ can be defined in a way
which is compatible with the general construction of Section
\ref{sec:basic}.

Choose any element $s_0 \in S$. Then $S_0: = \{s_0 s: s \in S\}$
is a generating set for $W^+$ which carries a simple
presentation~\cite[\S IV.1, Ex. 9]{Bou} and a Coxeter-like
structure~\cite{BRR}. Let us write
\[ \T_0 := \{ s_0 t:\ t\in \T\}. \]
Given a pair $(G, A)$ of a group $G$ and generating set $A$, we
say that an element $g \in G$ is an {\it odd palindrome} if there
is an $(A \cup A^{-1})^*$-word $(a_1,\ldots,a_\ell)$ for $g$ such
that $\ell$ is odd and $a_i = a_{\ell - i + 1}$ for every index
$i$. For example, the set of odd palindromes for $(W, S)$ is equal
to $\T$.

\begin{claim}
The set of odd palindromes for $(W^+, S_0)$ is equal to $\T_0$.
\end{claim}

\noindent
\emph{Proof.}
Let $w$ be an odd palindrome in $(W^+,S_0)$. Then $s_0 w$ is an
odd palindrome in $(W, S)$ and hence a reflection in $\T$.
Conversely, since $s_0$ is an involution, for every reflection $t=
s_{i_1} s_{i_2} s_{i_3} s_{i_4} s_{i_5} \cdots  s_{i_4} s_{i_3}
s_{i_2} s_{i_1}\in \T$,
$$
s_0 t= (s_0 s_{i_1}) (s_{i_2} s_0) (s_0 s_{i_3})  (s_{i_4} s_0)
(s_0 s_{i_5}) \cdots  (s_{i_4} s_0) (s_0 s_{i_3}) (s_{i_2} s_0)
(s_0 s_{i_1})
$$
is an odd palindrome in $(W^+, S_0)$.
\qed

\bigskip
Odd palindromes in alternating subgroups play a role which is
analogous to that played by reflections in Coxeter groups~\cite[\S
2.5, \S 3.5]{BRR}. This leads to the following definition of
absolute order to alternating subgroups.

\begin{defn} \rm
Given a simple reflection $s_0 \in S$, the {\it (left) absolute
order} $\leq_{\T_0}$ on the alternating subgroup $W^+$ of $W$ is
defined as the reflexive and transitive closure of the relation
consisting of the pairs $(u, v)$ of elements of $W^+$ for which
$\ell_{\T_0}(u) < \ell_{\T_0}(v)$ and $v = \tau u$ for some $\tau
\in \T_0$.
\end{defn}

The absolute order on $W^+$, which we will denote by $\Abs_0(W^+)$,
depends on the choice of $s_0$: non-conjugate simple reflections
determine non-isomorphic absolute orders on $W^+$. For example,
the absolute order on $B_3^+$ which is determined by the choice of
the adjacent transposition $s_0 = (1 \ 2) (-1 \, -2)$ is not
isomorphic to the one determined by the choice $s_0 = (1 \, -1)$.
However, the rank generating function is independent of the choice
of $s_0$. This  will be proved by considering the action of $W$ on
cosets of $\langle s_0 \rangle$, the subgroup generated by $s_0$.

Here are some basic lemmas on the absolute lengths $\ell_\T$ and
$\ell_{\T_0}$ which will be used in the proof. For $w \in W$ we set
$w^{s_0} := s_0 w s_0$.

\begin{lemma}\label{lemma1}
For every $w\in W^+$ we have
\[ \ell_{\T_0}(w) \ = \ \cases{
\ell_{\T}(w),  & if $\ell_{\T_0}(w)$ is even \cr
\ell_{\T}(w)-1, & if $\ell_{\T_0}(w)$ is odd. } \]
\end{lemma}

\noindent
\emph{Proof.}
Let $w = t_1 \cdots t_{\ell}$ be a $\T$-word for $w$ of length $\ell
:= \ell_\T(w)$. Since $w \in W^+$, the number $\ell$ is even and we
may write
\[ w \ = \ t_1 s_0 s_0 t_2 s_0 s_0 t_3\cdots t_{\ell-1}
s_0 s_0 t_{\ell} \ = \ s_0 t_1^{s_0} s_0 t_2 \cdots s_0
t_{\ell-1}^{s_0} s_0 t_{\ell}. \]
This proves that
\begin{equation}\label{alter-eq1}
\ell_{\T_0}(w) \, \le \, \ell \, = \, \ell_\T (w).
\end{equation}
Suppose that $\ell_{\T_0}(w)=2m$ is even. Then we may write
\[ w \ = \ s_0 t_1 \cdots s_0 t_{2m} \ = \ \prod \limits_{i=1}^{m}
\, t_{2i-1}^{s_0} t_{2i}, \]
with $t_i\in \T$ for each index $i$. Thus $\ell_\T (w) \le 2m =
\ell_{\T_0}(w)$ and the proof follows in this case. Finally, if
$\ell_{\T_0}(w) = 2m+1$ is odd, then we may write
\[ w \ = \ s_0 t_1 \cdots s_0 t_{2m+1} \ = \
s_0 t_1 \prod \limits_{i=1}^m \, t_{2i}^{s_0} \, t_{2i+1}, \]
with $t_i\in \T$ for each $i$. This shows that
\begin{equation} \label{eq3}
\ell_\T (w) \, \le \, 2m+2 \, = \, \ell_{\T_0}(w)+1.
\end{equation}
Combining (\ref{alter-eq1}) with (\ref{eq3}) yields  $\ell_{\T_0}(w)
\le \ell_\T (w)\le \ell_{\T_0}(w) + 1$. Since $\ell_{\T}(w)$ and
$\ell_{\T_0}(w)$ have distinct parities, we conclude that $\ell_\T
(w) = \ell_{\T_0}(w) + 1$ and the proof follows in this case too.
\qed

\begin{lemma} \label{lemma2}
For every $w\in W^+$, the following conditions are equivalent:
\begin{itemize}
\itemsep=0pt
\item[{\rm (i)}]
$\ell_{\T_0}(w)$ is even.
\item[{\rm (ii)}]
$\ell_{\T}(w) < \ell_\T (s_0 w)$.
\item[{\rm (iii)}]
$\ell_{\T}(w)<\ell_\T ( ws_0)$.
\end{itemize}
\end{lemma}

\noindent
\emph{Proof.}
Since the absolute length is invariant under conjugation
(Fact~\ref{absolute-length-properties} (e)), we have $\ell_\T (s_0 w)
= \ell_\T(w s_0)$ and hence it suffices to prove that ${\rm (i)}
\Leftrightarrow {\rm (ii)}$.

Suppose first that $\ell_{\T}(w)>\ell_\T (s_0 w)$. We note that $\ell_\T
(s_0 w)$ is an odd number, since $w\in W^+$, say $\ell_\T (s_0 w) =
2m+1$, and let $t_1 \cdots t_{2m+1}$ be a reduced $\T$-word for $s_0
w$. Then $s_0 t_1 \cdots t_{2m+1}$ is a reduced $\T$-word for $w$ and
$\ell_\T(w) = 2m+2$. Since
\[ w \ = \ s_0 t_1 \prod\limits_{i=1}^{m} \, (s_0 t_{2i}^{s_0}) (s_0
t_{2i+1}), \]
we have $\ell_{\T_0}(w) \le 2m+1$. On the other hand, we have
$\ell_{\T_0}(w) \ge \ell_\T(w) - 1 = 2m+1$ by Lemma~\ref{lemma1}. Thus
$\ell_{\T_0}(w) = 2m+1$ and, in particular, $\ell_{\T_0}(w)$ is odd.
This proves the implication ${\rm (i)} \Rightarrow {\rm (ii)}$.

Conversely, suppose that $\ell_{\T_0}(w)$ is odd. Then the proof of
Lemma~\ref{lemma1} shows that there is a reduced $\T$-word for $w$
which starts with $s_0$. This implies that $\ell_{\T}(w) > \ell_\T
(s_0 w)$ and hence ${\rm (ii)} \Rightarrow {\rm (i)}$. \qed

\bigskip
Let us denote by $\langle s_0\rangle$ the two-element subgroup of $W$
generated by $s_0$. We recall that the absolute length function on
$W / \langle s_0 \rangle$ is determined by Definition \ref{def:main2}.

\begin{corollary} \label{alternating-cor3}
We have $\ell_{\T}(w \langle s_0\rangle) = \ell_{\T_0}(w)$ for every
$w \in W^+$.
\end{corollary}

\noindent
\emph{Proof.}
By definition of $\ell_{\T}(w \langle s_0\rangle)$ and Lemma~\ref{lemma2}
we have
\[ \ell_{\T}(w \langle s_0\rangle) \ = \ \min \, \{\ell_\T(w),
\ell_\T(ws_0)\} \ = \ \cases{
\ell_{\T}(w), & if $\ell_{\T_0}(w)$ is even \cr
\ell_{\T}(w) - 1, & if $\ell_{\T_0}(w)$ is odd } \]
and the result follows from Lemma~\ref{lemma1}.
\qed

\begin{proposition}\label{iso}
The orders $\Abs_0(W^+)$ and $\Abs(W/\langle s_0\rangle)$ are
isomorphic.
\end{proposition}

\noindent
\emph{Proof.}
We consider the map $\varphi: W^+ \longrightarrow W / \langle
s_0 \rangle$ defined by
\[ \varphi(w) \ := \ \cases{
w  \langle s_0 \rangle,  & if $\ell_{\T_0}(w)$ is even \cr
w^{s_0} \langle s_0 \rangle, & if $\ell_{\T_0}(w)$ is odd} \]
for $w \in W^+$. We will show that $\varphi$ is the required
isomorphism of absolute orders.

We first note that conjugation by $s_0$ is an automorphism on
both $W$ and $W^+$ which preserves the lengths $\ell_\T$ and
$\ell_{\T_0}$, respectively. Corollary~\ref{alternating-cor3}
then implies that
\begin{equation}\label{varphi-length-invariance}
\ell_\T(\varphi(w)) \ = \ \ell_{\T_0}(w)
\end{equation}
for every $w\in W^+$. Since the map $\pi: W^+ \longrightarrow W /
\langle s_0 \rangle$ defined by $\pi(w) = w \langle s_0 \rangle$
is a bijection, we may conclude that $\varphi$ is a bijection as
well. Thus, it remains to show that the following conditions are
equivalent for $u, v \in W^+$:
\begin{itemize}
\itemsep=0pt
\item[(a)]
$u$ is covered by $v$ in $\Abs_0(W^+)$,
\item[(b)]
$\varphi (u)$ is covered by $\varphi(v)$ in $\Abs(W/\langle
s_0\rangle)$.
\end{itemize}
Using the definitions of the relevant absolute orders, we find that
\begin{itemize}
\itemsep=0pt
\item[(a)]
$\Leftrightarrow \, v = \tau u$ for some $\tau \in \T_0$ and
$\ell_{\T_0} (u) < \ell_{\T_0}(v)$
\item[]
$\Leftrightarrow \, v = s_0 t u$ for some $t \in \T$ and
$\ell_{\T_0} (u) < \ell_{\T_0}(v)$
\item[]
$\Leftrightarrow v^{s_0} \langle s_0 \rangle = t u \langle s_0
\rangle$ for some $t \in \T$ and $\ell_{\T_0} (u) < \ell_{\T_0}
(v)$
\item[]
$\Leftrightarrow v \langle s_0 \rangle = t^{s_0} u^{s_0} \langle
s_0 \rangle$ for some $t \in \T$ and $\ell_{\T_0} (u) < \ell_{\T_0}
(v)$
\end{itemize}
and
\begin{itemize}
\itemsep=0pt
\item[(b)]
$\Leftrightarrow \, \varphi(v) = t \varphi(u)$ for some $t \in \T$
and $\ell_{\T} (\varphi(u)) < \ell_{\T} (\varphi(v))$.
\end{itemize}
The claim that (a) $\Leftrightarrow$ (b) follows from the previous
equivalences, (\ref{varphi-length-invariance}) and the definition
of the map $\varphi$.
\qed

\bigskip
The following statement extends \cite[Theorem 7.2]{Rotbart} from
the case of symmetric groups to that of all finite Coxeter groups.

\begin{corollary}
For every finite Coxeter group $W$ we have
\begin{equation}\label{alternating-factor}
\sum\limits_{w\in W^+} q^{\ell_{\T_0}(w)} \ = \ {W_\T(q) \over 1+q}
\ = \ \prod\limits_{i=2}^d (1+e_i q),
\end{equation}
where $d$ is the Coxeter rank and $1 = e_1, e_2,\dots,e_d$ are
the exponents of $W$.
\end{corollary}

\noindent
\emph{Proof.}
Proposition~\ref{iso} implies that
\[ \sum\limits_{w\in W^+} q^{\ell_{\T_0}(w)} \ = \ (W/\langle
s_0\rangle)_\T(q). \]
Since $\langle s_0 \rangle$ is a modular subgroup of $W$ (see
Example~\ref{ex:modular1} (a)), we have
\[ (W/\langle s_0\rangle)_\T(q) \ = \ {W_\T(q)\over \langle s_0
\rangle_{\T}(q)} \ = \ {W_\T(q)\over 1+q} \]
by Proposition~\ref{prop:fact} and the first equality in
(\ref{alternating-factor}) follows. The second equality is a
restatement of (\ref{eq:rankgenW}).
\qed

\bigskip
Another description of $\Abs_0(W^+)$ can be given as follows. Let
us write $R_0 := \{ w \in W: \ell_\T(ws_0) > \ell_\T(w)\}$. The proof
of Proposition~\ref{prop:modular-ideal} shows that $R_0$ is an order
ideal of $\Abs(W)$.

\begin{corollary}\label{iso2}
The absolute order $\Abs_0(W^+)$ is isomorphic to $(R_0, \leq_\T)$.
\end{corollary}

\noindent
\emph{Proof.}
The proof of Proposition~\ref{prop:modular-ideal} shows that
$\Abs(W / \langle s_0\rangle)$ is isomorphic to $(R_0,\leq_\T)$.
The result follows from this statement and Proposition~\ref{iso}.
\qed

\section{Remarks on ordered tuples}
\label{sec:tuples}

This section briefly discusses the action of the symmetric group
$S_n$ on the set $X_{n, k}$ of ordered $k$-tuples of pairwise
distinct elements of $\{1, 2,\dots,n\}$, as well as a
generalization. The stabilizer $S_{n-k}$ of this action is a
modular reflection subgroup of $S_n$ (see Example~\ref{ex:SkBk}).
Therefore, by Proposition~\ref{prop:fact} we have
\[ {X_{n,k}}_\T(q) \ = \ {{S_n}_\T(q) \over
{S_{n-k}}_\T(q)} \ = \ \prod\limits_{i=n-k}^{n-1} (1+iq). \]

By a classical result of Hurwitz~\cite{Hurwitz} (see
also~\cite{Denes, Strehl}), there is a one-to-one correspondence
between the maximal chains of any maximal interval of $\Abs(S_n)$
and labeled trees of order $n$. The following generalization of
this statement on the enumeration of maximal chains of $X_{n, k}$
is possible. We will denote by $d_\Gamma (v)$ the valency (i.e.,
number of neighbors) of a node $v$ of a labeled tree $\Gamma$ of
order $n$.

\begin{proposition} \label{lemma-t1}
For all integers $1 \le k < n$, the number of maximal chains of
$\Abs(X_{n,k})$ is equal to
\[ k! \ \sum\limits_{\Gamma} \, (n-k)^{d_\Gamma (v_0)}, \]
where the sum runs over all trees $\Gamma$ on the node
set $\{v_0, v_1,\dots,v_k\}$.
\end{proposition}

The proof of this statement will be given elsewhere. The special
case $k=n-1$ is equivalent to Hurwitz's theorem.

\medskip
Combining Propositions~\ref{iso} and~\ref{lemma-t1}, we get the
following statement.

\begin{corollary}
The number of maximal chains of the absolute order on the alternating
group of $S_n$ is equal to
\[ (n-2)! \ \sum\limits_{\Gamma} \, 2^{d_\Gamma(v_0)}, \]
where the sum runs over all trees $\Gamma$ on the node set
$\{v_0, v_1,\dots,v_{n-2} \}$.
\end{corollary}

The previous setting has a natural extension to wreath product
actions on ordered colored tuples. Recall from~\cite{Shi1, Shi2}
the absolute order on the complex reflection group $G(r,n) = \ZZ_r
\wr S_n$; absolute length and order are naturally defined with
respect to the set $\T$, consisting of all elements
(pseudoreflections) of finite order fixing a hyperplane. Let
$X_{r,n,k} := \{(a_1,\dots,a_k): \forall i \ a_i \in \ZZ_r \times
\ZZ_n\}$ be the set of ordered $k$-tuples of letters in an alphabet
of size $n$ which are $r$-colored. Then $G(r,n)$ acts naturally on
$X_{r,n,k}$, with stabilizer $G(r,n-k)$. By extending
Propositions~\ref{prop:transitivity} and \ref{prop:maxmodular}, one
can prove that the subgroup $G(r,n-k)$ is a modular subgroup of
$G(r,n)$ for $1 \le k\le n$. Hence, by (the extension of)
Proposition~\ref{prop:fact} we have
\[ {X_{r,n,k}}_\T(q) \ = \ {{G(r,n)}_\T(q) \over
{G(r,n-k)}_\T(q)} \ = \ \prod\limits_{i=n-k}^{n-1} (1+riq). \]

\bigskip
\noindent
{\bf Acknowledgements.} The authors wish to thank the anonymous
referees for several helpful comments.


\begin{thebibliography}{99}

\bibitem{Armstrong}
D.~Armstrong,
{\it Generalized noncrossing partitions and combinatorics of
Coxeter groups}, Mem. Amer. Math. Soc. 202 (2009), no. 949.

\bibitem{AK}
C.A.~Athanasiadis and M.~Kallipoliti,
{\it The absolute order on the symmetric group, constructible
partially ordered sets and Cohen-Macaulay complexes},
J. Combin. Theory Series A {\bf~115} (2008), 1286-–1295.

\bibitem{BarceloIhrig98}
H.~Barcelo and E.~Ihrig,
{\it Modular elements in the lattice $L(A)$ when $A$ is a real
reflection arrangement},
Selected papers in honor of Adriano Garsia (Taormina, 1994),
Discrete Math. {\bf~193} (1998), 61–-68.

\bibitem{BarceloIhrig99}
H.~Barcelo and E.~Ihrig,
{\it Lattices of parabolic subgroups in connection with
hyperplane arrangements},
J. Algebraic Combin. {\bf~9} (1999), 5--24.

\bibitem{Ben-Ari}
R.~Ben-Ari,
{\it Combinatorial parameters on matchings in complete
graphs},
M. Sc. Thesis, Bar-Ilan University, 2010.

\bibitem{Be}
D.~Bessis,
{\it The dual braid monoid},
Ann. Sci. Ecole Norm. Sup. {\bf~36} (2003), 647--683.

\bibitem{BjB}
A.~Bj\"orner and F.~Brenti,
Combinatorics of Coxeter groups,
Graduate Texts in Mathematics {\bf~231},
Springer, New York, 2005.

\bibitem{Bou}
N.~Bourbaki,
Lie Groups and Lie Algebras,
Springer 2002.

\bibitem{BW}
T.~Brady and C.~Watt,
{\it $K(\pi,1)$'s for Artin groups of finite type},
in {\it Proceedings of the Conference on Geometric and Combinatorial
Group Theory}, Part I (Haifa, 2000),
Geom. Dedicata {\bf~94} (2002), 225–-250.

\bibitem{BRR}
F.~Brenti, V.~Reiner and Y.~Roichman,
{\it Alternating subgroups of Coxeter groups},
J. Combin. Theory Series A {\bf~115} (2008), 845–-877.

\bibitem{Denes}
J.~D\'enes,
{\em The representation of a permutation as the
product of a minimal number of transpositions and its connection
with the theory of graphs},
Publ.\ Math.\ Inst.\ Hungar.\ Acad.\ Sci.~{\bf 4} (1959), 63-–71.

\bibitem{Dyer01}
M.J.~Dyer,
{\it On minimal lengths of expressions of Coxeter group elements
as products of reflections},
Proc. Amer. Math. Soc. {\bf~129} (2001), 2591–-2595.

\bibitem{HHN}
C.~Hernando, F.~Hurtado and M.~Noy,
{\it Graphs of non-crossing perfect matchings},
Graphs Combin. {\bf~18} (2002), 517–-532.


\bibitem{Humphreys}
J.E.~Humphreys,
Reflection groups and Coxeter groups,
Cambridge Studies in Advanced Mathematics {\bf~29},
Cambridge University Press, Cambridge, 1990.

\bibitem{Hurwitz}
A.~Hurwitz,
{\it Ueber Riemann'sche Fl\"achen mit gegebenen Verzweigungspunkten}
(German), Math. Ann. {\bf~39} (1891), 1--60.

\bibitem{Ka}
M.~Kallipoliti,
{\it The absolute order on the hyperoctahedral group},
J. Algebraic Combin. {\bf~34} (2011), 183-–211.

\bibitem{Ku93}
J.P.S.~Kung,
{\it Sign-coherent identities for characteristic polynomials of
matroids},
Combin. Prob. Comput. {\bf~2} (1993), 33-–51.

\bibitem{OS80}
P.~Orlik and L.~Solomon,
{\it Unitary reflection groups and cohomology},
Invent. Math. {\bf~59} (1980), 77--94.

\bibitem{OT}
P.~Orlik, H.~Terao,
Arrangements of Hyperplanes,
Springer-Verlag, New York, 1992.

\bibitem{Rains-Vazirani}
E.M.~Rains and M.J.~Vazirani,
{\it Deformations of permutation representations of Coxeter groups},
{\tt arXiv:1008.1037}.

\bibitem{Rotbart}
A.~Rotbart,
{\it Generator Sets for the Alternating Group},
S\'{e}m. Lothar. Combin. {\bf~65} (2011),
Article B65b, 16pp (electronic).

\bibitem{Shi1}
J.-Y.~Shi,
{\it Formula for the reflection length of elements in the group
$G(m,p,n)$},
J. Algebra {\bf~316} (2007), 284–-296.

\bibitem{Shi2}
J.-Y.~Shi,
{\it Reflection ordering on the group $G(m,p,n)$},
J. Algebra {\bf~319} (2008), 4646–-4661.

\bibitem{Sloane}
N.J.A.~Sloane and S.~Plouffe,
The Encyclopedia of Integer Sequences, Academic Press, 1995.
Sequences: 
A007559, A007696.

\bibitem{Sta71}
R.P.~Stanley,
{\it Modular elements of geometric lattices},
Algebra Universalis {\bf~1} (1971), 214-–217.

\bibitem{StaEC1}
R.P.~Stanley,
Enumerative Combinatorics, vol.~1,
Wadsworth \& Brooks/Cole, Pacific Grove, CA, 1986;
second printing,
Cambridge Studies in Advanced Mathematics {\bf~49},
Cambridge University Press, Cambridge, 1997.

\bibitem{Stanley}
R.P.~Stanley,
{\it An Introduction to Hyperplane Arrangements},
in {\it Geometric Combinatorics}, 389–-496,
IAS/Park City Math. Series {\bf~13},
Amer. Math. Soc., Providence RI, 2007.

\bibitem{Strehl}
V.~Strehl,
{\it Minimal transitive products of transpositions -- the
reconstruction of a proof of A. Hurwitz},
S\'{e}m. Lothar. Combin. {\bf~37} (1996),
Article S37c, 12pp (electronic).

\end{thebibliography}
\end{document}